\documentstyle[12pt]{article}

\input amssym.def %
\input amssym.tex

\newtheorem{dfn}{Definition}[section]

\newtheorem{thm}[dfn]{Theorem}
\newtheorem{lem}[dfn]{Lemma}

\newtheorem{prop}[dfn]{Proposition}

\newtheorem{ex}[dfn]{Example}

\newtheorem{extprob}[dfn]{Extension Problem}

\def\proof{\par\medskip\noindent{\it Proof: }}

\def\ra{\rightarrow}

\def\R{{\Bbb R}}
\def\H{{\Bbb H}}
\def\Z{{\Bbb Z}}

\def\eps{\epsilon}

\def\ga{\gamma}

\def\Si{\Sigma}
\def\si{\sigma}

\def\D{\partial}

\def\smallprod{\mathop{\Pi}\limits}

\begin{document}
\title{3-manifolds with(out) metrics of nonpositive curvature}
\author{Bernhard Leeb \\ leeb@rhein.iam.uni-bonn.de}
\date{October 3, 1994}
\maketitle

\noindent
{\bf Abstract.}
In the context of Thurstons geometrisation program we address
the question which compact aspherical 3-manifolds admit
Riemannian metrics of nonpositive curvature.
We prove that a Haken manifold with, possibly empty, boundary of zero Euler
characteristic admits metrics of nonpositive curvature if the boundary is
non-empty or if at least
one atoroidal component occurs in its canonical topological decomposition.
Our arguments are based on Thurstons Hyperbolisation Theorem.
We give examples of closed graph-manifolds with linear gluing graph and
arbitrarily
many Seifert components which do not admit metrics of nonpositive curvature.

\section{Introduction}

It is known since the last century that a closed surface admits particularly
nice Riemannian metrics, namely metrics of constant curvature.
All aspherical surfaces, i.e.\ all surfaces besides
the sphere and the projective plane can be
given a metric of nonpositive curvature.

Thurston suggested a geometrisation procedure in dimension 3.
Here, constant curvature metrics are a too restricted class of model
geometries, and one allows more generally complete locally homogeneous
metrics.
There are 8 types of 3-dimensional geometries which can occur,
namely the model spaces are
$S^3,S^2\times R^1,R^3,Nil,Sol,\widetilde{SL(2,\R)},\H^2\times\R$ and $\H^3$.
The classification is due to Thurston \cite{Thurston1} and a detailed
exposition can be found in Scott's article \cite{Scott}.
Despite the large supply of model geometries
it is far from being true that any compact 3-dimensional manifold $M$ is
{\em geometric}, that is, admits a geometric structure.
One must cut $M$ into suitable pieces.
The {\em canonical topological decomposition} is obtained as follows:
According to Kneser, there is a maximal decomposition of $M$ as a
finite connected sum.
The summands are called prime manifolds and their homeomorphism type is
determined by $M$ if $M$ is orientable,
whereas the decomposition itself is in general not unique (Milnor).
The aspherical prime manifolds can be further decomposed by cutting along
incompressible embedded tori and Klein bottles, see
\cite{JacoShalen,Johannson}.
The components which one obtains in this second step are Seifert fibered or
atoroidal.
Thurstons {\em Geometrisation Conjecture} asserts that one can put a
geometric structure of unique type on each piece of the canonical
topological decomposition of $M$.
It is well-known that a closed 3-manifold is Seifert if and
only if it can be given one of the geometries different from $Sol$ and $\H^3$,
see \cite{Scott}.
The atoroidal components conjecturally admit a hyperbolic structure
and this has been proven by Thurston \cite{Thurston2} for components which are
Haken, i.e.\ contain a closed incompressible surface.
In particular, all Haken 3-manifolds are geometrisable.

We are interested in those 3-manifolds $M$ which have a chance
to admit metrics of nonpositive curvature.
In the context of Thurstons geometrisation program,
we address in this paper the

\medskip\noindent
{\bf Question:}
{\em Which compact aspherical 3-manifolds $M$
admit Riemannian metrics of nonpositive sectional curvature?}

\medskip\noindent
The situation is well-understood for geometric 3-manifolds.
If $M$ is Haken and non-geometric, then each component
in the canonical decomposition admits a geometric structure which is modelled
on one of the nonpositively curved geometries $\R^3,\H^2\times \R$ and $\H^3$.
The known obstructions to the existence of a
nonpositively curved metric on $M$ vanish,
i.e.\ all solvable subgroups of $\pi_1(M)$ are virtually abelian and
centralisers virtually split.
It is therefore only natural to ask whether one can put compatible
nonpositively curved metrics on all components of $M$ simultaneously.
We show that, as suggested by the geometrisation program, non-geometric Haken
manifolds indeed generically admit metrics of nonpositive curvature.
More precisely, we prove the following existence results:

\medskip\noindent
{\bf Theorem \ref{existence: graph}}
{\em
Suppose that $M$ is a graph-manifold, i.e.\ contains only Seifert components,
and has non-empty boundary.
Then there exists a Riemannian metric of nonpositive curvature on $M$.}

\medskip\noindent
And, relying on Thurstons Hyperbolisation Theorem \cite{Thurston2}:

\medskip\noindent
{\bf Theorem \ref{with atoroidal piece}}
{\em
Let $M$ be a Haken manifold with, possibly empty, boundary of zero Euler
characteristic.
Suppose that at least one atoroidal component occurs in its canonical
decomposition.
Then $M$ admits a Riemannian metric of nonpositive curvature.}

\medskip\noindent
The Riemannian metrics can be chosen smooth and flat near the boundary.
Our existence results yield new examples of closed nonpositively curved
manifolds of geometric rank one with non-hyperbolic fundamental group.

Differently from the situation in dimension 2, not all closed
aspherical 3-manifolds admit metrics of nonpositive curvature.
Examples can already be found among geometric manifolds,
namely quotients of $Nil$, $Sol$ or $\widetilde{SL(2,\R)}$.
We show that also non-geometric Haken manifolds cannot always be
equipped with a metric of nonpositive curvature.
We give the following

\medskip\noindent
{\bf Example \ref{non-existence}}
{\em
There are closed graph-manifolds glued from arbitrarily many Seifert
components which do {\bf not} admit metrics of nonpositive curvature.}

\medskip\noindent
Note that nevertheless
the fundamental groups of such manifolds are
nonpositively curved {\em in the large}:
We prove in \cite{graph} that the fundamental group of a closed Haken
3-manifold
is quasi-isometric to $Nil$, $Sol$ or the fundamental group of a closed
nonpositively curved 3-manifold.

\medskip
The paper is organised as follows:
In section \ref{geometric components},
we explain that any nonpositively curved
metric on a Haken manifold arises from nonpositively curved metrics on its
geometric components (section \ref{reduction})
and is {\em rigid} on Seifert components.
The atoroidal components, on the other hand, are {\em flexible}
(section \ref{atoroidal pieces}).
The extent of rigidity for the Seifert pieces is described in section
\ref{Seifert}.
The rigidity is responsible for the non-existence examples
(section \ref{example})
but leaves sufficient degrees of freedom
(Proposition \ref{flexibility of Seifert})
to produce nonpositively curved
metrics on Haken manifolds with atoroidal components {\em or} non-empty
boundary (section \ref{Existence results}).

We refer to \cite{CheegerEbin} for an introduction to the geometry of
nonpositive curvature
and to \cite{Jaco} for concepts from 3-dimensional topology.

\medskip\noindent
{\bf Acknowledgements.}
The results in this paper are part of my thesis \cite{thesis}.
I benefitted from the inspiring environment provided by the geometry
group at the University of Maryland.
I am particularly grateful to
Bill Goldman, Karsten Grove and John Millson for their advice and
encouragement during my time as a graduate student and to
Misha Kapovich for many helpful discussions.

\section{Nonpositively curved metrics on the geometric components}
\label{geometric components}

\subsection{Reduction to geometric components}
\label{reduction}

In this section $M$ will always denote a compact Haken 3-manifold with
boundary of zero Euler characteristic.
We recall that a compact smooth aspherical 3-manifold is {\em Haken} if
it contains a closed incompressible surface,
that is,
a closed smooth embedded 2-sided surface whose fundamental group is
infinite and injects via the canonical inclusion homomorphism,
cf.\ \cite{Jaco}.

According to \cite{JacoShalen,Johannson},
there is a canonical topological decomposition of $M$ into compact 3-manifolds
which are Seifert fibered or atoroidal.
It is obtained by cutting $M$ along finitely many
disjoint closed embedded incompressible tori and Klein bottles $\Si_i$
which we call the {\em splitting} or {\em decomposing} surfaces.
We refer to the pieces of the decomposition as {\em geometric components}
because they
can be equipped with canonical geometric structures.
If the topological decomposition of $M$ is non-trivial,
the geometric components have non-empty incompressible boundary.
The Seifert pieces then admit a $\R^3$- or $\H^2\times\R$-structure
(compare \cite{Scott})
whereas the atoroidal pieces admit $\H^3$-structures \cite{Thurston2}.
A minimal topological decomposition of $M$ is unique up to isotopy.

Assume now that $M$
carries a Riemannian metric $g$ of nonpositive sectional curvature.
We only allow metrics where the boundary is totally-geodesic and (hence) flat.
In this geometric situation,
there is an analoguous geometric decomposition of $(M,g)$ along
totally-geodesi\-cal\-ly embedded flat surfaces,
cf.\ \cite{thesis,deco}.
Since a minimal topological decomposition of $M$ is unique,
it can therefore be geometrised in the presence of $g$:

\begin{lem}[\cite{thesis,deco}]
There is an isotopy of $M$ which moves
the splitting surface $\Si:=\cup_i\Si_i$ to a totally-geodesically
embedded flat surface in $M$.
\end{lem}

It follows that each metric of nonpositive curvature on $M$ arises
modulo isotopy in the following way:
Put suitable flat metrics on the decomposing surfaces $\Si_i$ and
extend them to nonpositively curved metrics on the geometric components of the
topological decomposition.
It is therefore crucial to understand the following question:

\begin{extprob}
\label{extension problem}
Let $X$ be a geometric component.
Which flat metrics on the boundary $\D X$ can be extended to metrics of
nonpositive curvature on $X$?
\end{extprob}

Atoroidal and Seifert components behave differently with respect to the
extendability of metrics.
As discussed below,
atoroidal components are so {\em flexible} that the extension problem can be
solved for all flat metrics prescribed on the boundary,
whereas the rigidity of nonpositively curved Seifert manifolds puts
restrictions on the solvability of the extension problem for Seifert
components.
The interdependence of flat metrics on the boundary tori of a nonpositively
curved Seifert manifold is the source of obstructions to the existence of
nonpositively curved metrics on graph-manifolds.

\subsection{Nonpositively curved metrics on atoroidal manifolds}
\label{atoroidal pieces}

Let $X$ be an atoroidal component of $M$.
We already mentioned a deep result of Thurston \cite{Thurston2},
his Hyperbolisation Theorem,
that the interior $int X$ of
$X$ admits a complete metric $g_0$ of constant negative curvature, say $-1$,
and finite volume.
The hyperbolic structure $g_0$ is in fact unique according to Mostows
rigidity theorem.
By modifying $g_0$ we can solve the extension problem:

\begin{prop}[Flexibility of atoroidal components]
\label{atoroidal comp}
Any flat metric on $\D X$ can be extended to a smooth nonpositively
curved metric on $X$ which is flat near the boundary.
\end{prop}

\proof
Let $h$ be a prescribed flat metric on $\D X$.
We will change the hyperbolic structure $g_0$ in order to extend $h$ to
all of $X$.
The ends of $(int X,g_0)$ are hyperbolic cusps:
Outside a suitable compact subset, the metric $g_0$ is isometric to a warped
product metric
\[ e^{-2t}g_{\D X}+dt^2 \]
on $\D X\times\R^+$ where $g_{\D X}$ is a flat metric on $\D X$.
The conformal type of $g_{\D X}$ is determined by $g_0$.
There is no relation between the metrics $g_{\D X}$ and $h$.
As a first step,
we adjust the conformal type of the cusps of $g_0$ to the prescribed metric
$h$.
To do so it suffices to allow more general warped product metrics:
after isotoping $h$ if necessary
we can diagonalise $h$ with respect to $g_{\D X}$,
that is, we write
$g_{\D X}$ and $h$ as
\[ g_{\D X}=dx^2+dy^2
\qquad\hbox{ and }\qquad h=a^2\cdot dx^2+ b^2\cdot dy^2 \]
where $dx^2$ and $dy^2$ are positive-semidefinite sub-Riemannian metrics
parallel with respect to $g_{\D X}$
with orthogonal one-dimensional kernels
and $a,b$ are positive functions constant on each boundary component
$\Si_i$.
To interpolate between the conformal types of $g_{\D X}$ and $h$, we put a
metric of the form
\[
e^{-2t}\bigl[\phi+(1-\phi)\cdot a\bigr]^2\cdot dx^2 +
e^{-2t}\bigl[\phi+(1-\phi)\cdot b\bigr]^2\cdot dy^2 +dt^2 \]
on $\D X\times\R^+$
where the smooth function $\phi:\R^+\to[0,1]$ is required to be
equal
to 1 in a neighborhood of 0 and equal to 0 in a neighborhood of $\infty$.
A curvature calculation shows that the sectional curvature of a
metric of this type is negative if the component functions
$ e^{-t}\bigl[\phi+(1-\phi)\cdot a\bigr]$ and
$ e^{-t}\bigl[\phi+(1-\phi)\cdot b\bigr]$
are strictly monotonically decreasing and convex.
This holds for all functions $\phi$ whose
first and second derivatives are
bounded by sufficiently small constants depending on $a$ and $b$.
Hence we can find a complete negatively curved metric $g_1$ on
$int X$ which is outside a suitable compact subset hyperbolic and isometric
to the warped product metric
\[ e^{-2t}h+dt^2 \]
on $\D X\times(T_0,\infty)$ for some $T_0\in\R$.

In the second step,
we replace $e^{-2t}$ by a convex and monotonically decreasing function
$\psi:(T_0,\infty)\to(0,\infty)$
which coincides with $e^{-2t}$ in a neighborhood of $T_0$ and is constant in a
neighborhood of $\infty$.
The curvature of the resulting complete metric $g_2$
is nonpositive because $\psi$ is convex.
After rescaling,
$(int X,g_2)$ is outside some compact subset isometric to a Euclidean cylinder
with base $(\D X,h)$.
We cut off the ends along cross sections of the cylinders and obtain a smooth
nonpositively curved metric on $X$ which is flat near the boundary and extends
$h$.
$\Box$

\subsection{Nonpositively curved metrics on Seifert manifolds}
\label{Seifert}

We address the Extension Problem \ref{extension problem} for Seifert
manifolds.
Let $X$ be a Seifert fibered manifold with non-empty incompressible boundary
and denote by $O$ its base orbifold.
We have the exact sequence
\[ 1\to\langle f\rangle\to\pi_1(X)\to\pi_1(O)\to1 \]
where the element $f$ is represented by a generic Seifert fiber.
We recall that
due to the cyclic normal subgroup of the fundamental group,
any nonpositively curved metric $g$ on $X$ must have a special form:
The universal cover decomposes as a Riemannian product
\begin{equation}
\label{1}
 (\tilde X,g)\cong \R\times Y
\end{equation}
where $Y$ is nonpositively curved, cf. \cite{Eberlein}.
The lines $\R\times\{y\}$ are the axes of the deck-transformation $f$,
and they project down to a Seifert fibration of $X$ by closed geodesics
of, apart from finitely many singular fibers, equal length.
The deck-action of $\pi_1(X)$ preserves the splitting (\ref{1}) and hence
decomposes as the product of a representation
\begin{equation}
\label{2}
\phi: \pi_1(X)\to Isom(\R)
\cong\R\rtimes\Z/2\Z
\end{equation}
and a cocompact discrete action of $\pi_1(O)$ on $Y$
which induces an orbifold-metric of nonpositive curvature on $O$.
Since $f$ acts on its axes by translations, the representation $\phi$
satisfies the condition
\begin{equation}
\label{3}
 \phi(f)\in\R\setminus\{0\} .
\end{equation}
We see that a nonpositively curved metric on $X$ correponds to the choice
of a nonpositively curved metric on $O$ and a representation (\ref{2})
satisfying (\ref{3}).
It is well-known that
the Seifert manifold $X$ does admit metrics of nonpositive curvature.
This follows, for instance, from Theorem 5.3 in \cite{Scott};
namely, the closed Seifert manifold obtained from doubling $X$ admits a metric
of nonpositive curvature which yields a representation (\ref{2}) for
$\pi_1(X)$ so that (\ref{3}) holds.
We call a Seifert manifold of {\em Euclidean} respectively {\em hyperbolic}
type according to whether its base orbifold $O$ admits flat or hyperbolic
metrics.
All nonpositively curved metrics on a Seifert manifold of Euclidean type are
flat.

Suppose in the remainder of this section that $X$ has hyperbolic type.
We want to determine which flat metrics $h$ on $\D X$ can be extended to a
nonpositively curved metric on $X$.
One necessary condition is immediate from the above discussion:
The closed geodesics in $(\D X,h)$ which are (homotopic to) Seifert fibers
must have equal lenghts.
This is sufficient in the non-orientable case:

\begin{lem}
\label{non-orientable Seifert}
If $X$ is non-orientable,
then any flat metric on $\D X$ so that Seifert fibers have equal length can be
extended to a nonpositively curved metric on $X$.
\end{lem}

In the orientable case,
there is a second necessary condition because the values of
the representation $\phi$ on the
boundary tori $T_i$ are interrelated.
To describe it we need some notation:
The flat metric $h$ on $\D X$ induces scalar products $\si_i$
on the abelian groups $\pi_1(T_i)\cong H_1(T_i,\Z)$.
The values of $\phi$ on $H_1(T_i,\Z)$ lie in $\R$
and they are related to $\si_i$ by the formula
\[ \si_i(f,\cdot)=\phi(f)\cdot\phi \qquad\hbox{ on $H_1(T_i,\Z)$.} \]
We choose in each group $H_1(T_i,\Z)$ an element $b_i$ which forms
together with $f$ a basis compatible with the orientation of $T_i$.
The compatibility condition for the flat metrics on the boundary tori is given
by

\begin{lem}
[Rigidity of Seifert components]
\label{orientable Seifert}
If $X$ is oriented,
there is a rational number $c$ which depends on the types of
singular fibers
and the choice of the basis elements $b_i$, so that the following is true:
A flat metric $h$ on $\D X$ can be extended to a nonpositively curved metric
on $X$ if and only if the Seifert fibers of $(\D X,h)$ have equal length and
the relation
\begin{equation}
\label{compatibility cond}
 \sum_i \si_i(f,b_i)=c\cdot\|f\|^2
\end{equation}
holds.
\end{lem}

Note that the nonpositively curved metrics on $X$ can always be chosen flat
near the boundary by adjusting the metric on the base-orbifold $O$.

The proofs are straight-forward.
By varying the hyperbolic structure on $O$ one can realize
arbitrary lengths for the boundary curves.
The possible values for the representation $\phi$ on the boundary $\D X$ can
be analysed as in \cite{Scott} by using a canonical presentation for
$\pi_1(X)$.
A proof for the non-orientable case is given in \cite{action}.

Observe that by varying the metric on the base orbifold $O$, exactly those
modifications of the induced flat metric on $\D X$ are possible which preserve
the scalar products
$\si_i(f,\cdot)$ with the fiber direction and, consequently,
the length $\|f\|$ of the fiber.
We will refer to
this way of changing the flat metric on $\D X$ as
{\em rescaling orthogonally to the Seifert fiber}.

In our construction of nonpositively curved metrics on graph-manifolds in
section \ref{Existence results} we need to extend flat metrics which are
prescribed only on part of the boundary.
Lemmata \ref{non-orientable Seifert} and \ref{orientable Seifert} imply that
the Extension Problem \ref{extension problem} with loose ends is always
solvable:

\begin{prop}
[Restricted flexibility of Seifert components]
\label{flexibility of Seifert}
Sup-
\par\noindent
pose that $X$ is a Seifert manifold with non-empty incompressible
boundary and hyperbolic base orbifold.
Let $h$ be a flat metric which is prescribed on some but not all of the
boundary components of $X$ so that Seifert fibers have equal lengths.
Then $h$ can be extended to a nonpositively curved metric on $X$ which is
flat near the boundary.
\end{prop}

For the sake of completeness, we show how one can deduce the proposition
from the existence of a single nonpositively curved metric $g$ on $X$.
This is done by modifying the representation $\phi:\pi_1(X)\to Isom(\R)$.
There is a simple closed loop in $O$ which separates the boundary from the
singularities and orientation reversing loops.
Accordingly $X$ can be obtained from gluing a Seifert manifold $Y$ and a
circle bundle $Z$ over a punctured sphere along one boundary component.
Denote by $a_0,\dots,a_k$ the compatibly oriented boundary curves of a section
of $Z$ so that $a_0$ corresponds to the boundary surface along which $Z$ is
glued to $Y$.
The fundamental group of $Z$ is presented by
\[ \Bigl\langle
f,a_0,\dots,a_k \Big| a_ifa_i^{-1}=f^{\eps_i},\smallprod_i a_i=1
\Bigr\rangle \]
with $\eps_i\in\{\pm1\}$ and $\Pi_i\eps_i=1$.
A new presentation $\phi'$ is obtained from $\phi$ as follows.
We leave the presentation $\phi$ unchanged on $\pi_1(Y)$ and choose $\phi'$ on
$Z$ so that $\phi'(a_i)$ is in the same component of $Isom(\R)$ as $\phi(a_i)$
and so that the conditions
\[ \phi'(f)=\phi(f),\qquad
\phi'(a_0)=\phi(a_0) \qquad\hbox{ and }\qquad
\smallprod_i\phi'(a_i)=\smallprod_i\phi(a_i)\]
hold.
This shows that we have the freedom to prescribe flat metrics on $\D X$ as
claimed in Proposition \ref{flexibility of Seifert}.

\section{Existence results}
\label{Existence results}

\subsection{Graph-manifolds with boundary}

There is no unanimous notion of graph-manifold available in the literature.
We choose the

\begin{dfn}
A {\bf graph-manifold} is a Haken manifold which contains only Seifert
components in its topological decomposition.
\end{dfn}

Any metric of nonpositive curvature on a graph-manifold $M$ is rigid in the
sense that it splits almost everywhere locally as a product,
see section \ref{Seifert}.
Nevertheless, the restricted flexibility of Seifert manifolds as stated in
Proposition \ref{flexibility of Seifert} suffices to construct
nonpositively curved metrics on $M$ in the presence of boundary:

\begin{thm}
\label{existence: graph}
Suppose that $M$ is a graph-manifold with non-empty boundary.
Then there exists a Riemannian metric of nonpositive curvature on $M$.
\par\noindent
{\bf Addendum.}
Moreover, let $\{\ga_i\}$ be a collection of homotopically non-trivial
simple loops
in the boundary, one on each component, so that none of them represents the
fiber of the adjacent Seifert component.
Then, given positive numbers $l_i$, the nonpositively curved metric on $M$ can
be chosen so that the loops $\ga_i$ are geodesics of length $l_i$.
\end{thm}

We will give examples of closed graph-manifolds without metrics of nonpositive
curvature in section \ref{Non-existence examples}.

\proof
The Seifert components of $M$ have non-empty incompressible boundary and hence
admit nonpositively curved metrics.
We can assume that each Seifert component of Euclidean type has one end, since
otherwise it would be homeomorphic to a torus or Klein bottle cross the unit
interval and hence be obsolete in the topological decomposition.
Furthermore, we may assume that no gluing of ends of Seifert components
identifies the Seifert fibers.

Let $n$ be the number of Seifert components of $M$ and suppose that the claim
has been proven for graph-manifolds with less than $n$ components.
Choose a Seifert component $X$ which contributes to the boundary of $M$.
It has hyperbolic type.
The complement $C$ of $X$ in $M$ consists of some Seifert manifolds $X_i^{eu}$
of Euclidean type with one end and a graph-manifold $M_0$ whose Seifert
components adjacent to $X$ have hyperbolic base.
We pick flat metrics on the Euclidean pieces $X_i^{eu}$ and, using the
induction assumption, a nonpositively curved metric on $M_0$.
We can arrange this metric $g_0$ on $C$ so that all closed
geodesics in the boundary which
are identified with fibers in $\D X$ via the gluing, have the same length $l$.
This is achieved on $M_0$ by rescaling orthogonally to the Seifert fiber,
as described in section \ref{Seifert}.
The metric $g_0$ determines via the gluing map a flat metric on $\D X\cap C$
with the property that Seifert fibers have equal length $l$.
It may also happen that boundary surfaces of $X$ are glued to each other.
We pick flat metrics on those which are compatible with the
identification and give length $l$ to the Seifert fibers as well.
Our choices prescribe a flat metric on a proper part of $\D X$ so that Seifert
fibers have equal length.
According to Proposition \ref{flexibility of Seifert},
this flat metric can be extended to a nonpositively curved metric on $X$.
Combined with $g_0$, this yields a smooth nonpositively curved metric on $M$.
The Addendum follows by rescaling orthogonally to the fiber and ordinary
rescaling.
$\Box$

\subsection{Haken manifolds with atoroidal components}

By now we have collected all ingredients for our main existence result.
We saw in
Proposition \ref{atoroidal comp}
that atoroidal manifolds are completely flexible,
and
the existence theorem \ref{existence: graph} for graph-manifolds with
boundary allows to put
nonpositively curved metrics on the complement of the atoroidal components of
a Haken manifold.
Combining these facts we get:

\begin{thm}[Existence in presence of an atoroidal piece]
\label{with atoroidal piece}
Let $M$ be a Haken manifold with, possibly empty, boundary of zero Euler
characteristic.
Suppose that at least one atoroidal component occurs in its canonical
decomposition.
Then $M$ admits a Riemannian metric of nonpositive curvature.
\end{thm}

We obtain new examples of closed nonpositively curved
manifolds which are not locally symmetric and do not admit metrics of
strictly negative curvature.
These manifolds have geometric rank one in the sense of Ballmann, Brin and
Eberlein.
The first 3-manifold examples of this kind were given by Heintze
\cite{Heintze} and Gromov \cite{Gromov}.

\section{Closed graph-manifolds}
\label{Non-existence examples}

In section \ref{Existence results},
we proved existence of nonpositively curved metrics on
Haken manifolds which have non-empty boundary or contain an atoroidal
component in their topological decomposition.
This leaves open the case of closed graph-manifolds which we address now.

Consider a closed non-geometric graph-manifold $M$.
As explained in section \ref{geometric components},
a metric $g$ of nonpositive curvature on $M$ induces flat metrics $g_i$ on the
decomposing surfaces $\Si_i$
and $g$ locally splits as a product on each Seifert component.
This rigidity interrelates the metrics $g_i$ and serves as the source of
obstructions to the existence of nonpositively curved metrics.
A nonpositively curved metric exists on $M$
if and only if we can put flat metrics on the Euclidean Seifert components and
the surfaces $\Si_i$ so that all orientable Seifert components with hyperbolic
base satisfy the linear compatibility condition stated in
Lemma \ref{orientable Seifert}.
This reduces the existence question to a purely algebraic problem which is,
however, in general fairly complicated.
In section \ref{example} below
we discuss a class of graph-manifolds with linear gluing graph
where the compatibility conditions
can be translated into a simple combinatorial criterion for finite
point configurations in hyperbolic plane.
Although a majority of these manifolds admit metrics of nonpositive
curvature, we will also obtain examples of non-existence.

The method in section \ref{example}
can be extended to the more general case when the dual graph to
the canonical decomposition is a tree (joint work with Misha Kapovich).
Examples of graph-manifolds with no metrics of nonpositive curvature
when the dual graph contains cycles are given in \cite{action};
they arise as mapping tori of reducible surface diffeomorphisms.
Recently, Buyalo and Kobelski \cite{Buyalo} announced the construction of a
numerical
invariant for graph-manifolds which detects the existence of a nonpositively
curved metric.

\subsection{An example}
\label{example}

We consider closed graph-manifolds $M$ built from two Seifert components
$X_0$ and $X_{n+1}$, each with one boundary torus
(denoted by $\partial_+X_0$ respectively $\partial_-X_{n+1}$),
and $n\geq0$ Seifert components $X_1,\dots,X_n$ with two boundary tori
(denoted by $\partial_{\pm}X_i$).
The $X_i$ are assumed to be trivial circle bundles
over compact orientable surfaces $\Sigma_i$
of genus $\geq1$ and
hence they are of hyperbolic type.
The graph-manifold $M$ is obtained by performing gluing homeomorphisms
\[ \alpha_i:\partial_+X_i\ra\partial_-X_{i+1} \qquad (i=0,\dots,n). \]
Denote by
$T_i:=\partial_+X_i=\partial_-X_{i+1}$
the decomposing tori of $M$.

We define an abstract free abelian group $A$ of rank 2
by identifying the homology groups $H_1(T_i,\Z)$ with each other in a
canonical way:
For $i=1,\dots,n$, we define the elements
$x_{i-1}\in H_1(T_{i-1},\Z)$ and $x_i\in H_1(T_i,\Z)$ to be equivalent iff
their images under the inclusion monomorphisms
\[H_1(\D_{\pm}X_i,\Z)\hookrightarrow H_1(X_i,\Z)\]
coincide.
The Seifert fibers of the components $X_i$ yield distinguished elements
$f_0,\dots,f_{n+1}\in A$.
In addition there are elements $b_0,b_{n+1}\in A$ which correspond to the
well-defined horizontal directions in the boundary tori of the one-ended
components; they can be defined as generators of the kernels of the natural
homomorphisms
\[ H_1(\D_+X_0,\Z)\to H_1(X_0,\Z) \] 
and
\[ H_1(\D_-X_{n+1},\Z)\to H_1(X_{n+1},\Z) .\]
Note that there are no well-defined horizontal directions in the boundary tori
of components with more than one end.
The collection of vectors $b_0,b_{n+1},f_0,$ $\dots,$ $f_{n+1}$
encodes the gluing maps and hence the topology of the closed graph-manifold
$M$.

To put flat metrics on the tori $T_i$ is equivalent to choosing scalar
products, i.e.\ positive-definite symmetric bilinear forms $\si_i$ on the
group $A$.
We recall that the space ${\cal H}$ of projective equivalence classes of
scalar products on $A$ carries a natural metric which is isometric to
hyperbolic plane.
A geodesic $c$ in ${\cal H}$ corresponds to a splitting
$V=L_1\oplus L_2$ of the 2-dimensional vector space
$V:=A\otimes_{\Z}\R$ into 1-dimensional subspaces;
namely $c$ consists of all (classes of) scalar products with respect to which
the splitting is orthogonal.
The ideal boundary points in the geometric compactification of ${\cal H}$
correspond to projective equivalence classes of positive-semidefinite forms
or, in other words,
to 1-dimensional subspaces of $V$.
The geometric boundary $\D_{geo}{\cal H}$ can hence be canonically
identified with the projective line ${\Bbb P}V$.
In particular, each element $a\in A$ can be interpreted as a point
$[a]\in\D_{geo}{\cal H}$.

As in Lemma \ref{orientable Seifert} we obtain the following necessary and
sufficient compatibility
conditions for the scalar products $\si_0,\dots,\si_n$ on $A$ so that the
corresponding flat metrics on the decomposing tori $T_i$ can be extended to a
nonpositively curved metric on $M$:
For the one-ended components we have
\[ \si_0(f_0,b_0)=0 \qquad\hbox{ and }\qquad \si_n(f_{n+1},b_{n+1})=0, \]
because the horizontal boundary curves representing $b_0,b_{n+1}$ are
homologically trivial in $X_0$ respectively $X_{n+1}$.
For the components with two ends we get
\[ \si_{i-1}(f_i,\cdot)=\si_i(f_i,\cdot) \qquad\hbox{ for $i=1,\dots,n$.} \]
This translates into geometric conditions relating
the ideal points $[b_0],[b_{n+1}],$ $[f_0],\dots,$ $[f_{n+1}]\in\D_{geo}{\cal
H}$ and the points $[\si_0],\dots,[\si_n]\in{\cal H}$:
\begin{enumerate}
\item
$[\si_0]$ lies on the geodesic with ideal endpoints $[b_0]$ and $[f_0]$.
$[\si_n]$ lies on the geodesic with endpoints $[b_{n+1}]$ and $[f_{n+1}]$.
\item
$[\si_{i-1}]$ and $[\si_i]$ lie on a geodesic asymptotic to $[f_i]$
($i=1,\dots,n$).
\end{enumerate}

\noindent
Note that this is a condition on the conformal types of the metrics on the
tori $T_i$.
But since the gluing graph does not contain cycles,
a collection of compatible conformal structures gives rise to a collection of
compatible flat metrics.
Our discussion yields the following

\medskip\noindent
{\bf Criterion.}
{\em
There exists a metric of nonpositive curvature on $M$
if and only if there is a configuration of points $[\si_0],\dots,[\si_n]$ in
${\cal H}$ which satisfies the above conditions 1 and 2.}

\medskip\noindent
Let us first consider the simplest case when $M$ is obtained from gluing two
Seifert pieces with one end:

\begin{ex}
If $n=0$, then $M$ admits a metric of nonpositive curvature if and only if
the gluing map $\alpha_0:\D_+X_0\to\D_-X_1$ preserves the canonical bases
of the first homology groups,
i.e.\ if
$\alpha_{0,\ast}:H_1(\D_+X_0,\Z)\to H_1(\D_-X_1,\Z)$ satisfies:
\[ \alpha_{0,\ast} \{\pm b_0,\pm f_0\} = \{\pm b_1,\pm f_1\} \]
\end{ex}

\proof
If a nonpositively curved metric exists on $M$,
then the splittings
$A=\langle b_0\rangle\oplus\langle f_0\rangle$ and
$A=\langle b_1\rangle\oplus\langle f_1\rangle$
must be orthogonal with respect to $\si_0$.
The sets
$\{\pm b_0,\pm f_0\}$ and $\{\pm b_1,\pm f_1\}$ contain the shortest
lattice vectors with respect to $\si_0$ and therefore coincide.
$\Box$

One can as well produce examples
with arbitrary number of Seifert components:
Let $C_0$ be the geodesic in ${\cal H}$ with endpoints $[f_0]$ and $[b_0]$.
For $i=1,\dots,n$ define $C_i$ to be the
union of all geodesics in ${\cal H}$ which are asymptotic to
$[f_i]$ and intersect $C_{i-1}$.
This is an increasing sequence of convex subsets of ${\cal H}$,
and $C_i$ consists of all conformal structures on $T_i$ which can be induced
by a nonpositively curved metric on $X_0\cup\dots\cup X_i$.
There exists a nonpositively curved metric on $M$ if and only if the geodesic
connecting $[f_{n+1}]$ and $[b_{n+1}]$ intersects $C_n$.
The contrary can easily be arranged,
choose for instance:
\begin{eqnarray*}
f_i &:=& f_0+i\cdot b_0 \qquad\hbox{ for $i=1,\dots,n+1$} \\
b_{n+1} &:=& f_0+(n+2)\cdot b_0
\end{eqnarray*}
Then $b_0,f_0,\dots,f_{n+1},b_{n+1}$ lie in this order on the
circle $\D_{geo}{\cal H}$ and no metric of
nonpositive curvature exists on the corresponding graph-manifold $M$.
This yields:

\begin{ex}
\label{non-existence}
There exist closed graph-manifolds with arbitrarily many Seifert components
which do {\bf not} admit metrics of nonpositive curvature.
\end{ex}


\begin{thebibliography}{Ch-E}
\addcontentsline{toc}{chapter}{Bibliography}

\bibitem[B-K]{Buyalo}
S.\ Buyalo and V.\ Kobelski,
{\em Geometrisation of graph-manifolds II: isometric states},
in preparation.

\bibitem[Ch-E]{CheegerEbin}
J.\ Cheeger and D.\ Ebin,
{\em Comparison theorems in Riemannian geometry},
North Holland 1975.

\bibitem[E]{Eberlein}
P.\ Eberlein,
{\em A canonical form for compact nonpositively curved manifolds
whose fundamental groups have nontrivial center},
Math.\ Ann.\ 260, (1982), vol.\ 1, 23-29.

\bibitem[G]{Gromov}
M.\ Gromov,
{\em Manifolds of negative curvature},
J.\ Diff.\ Geom.\ 13 (1978), 223-230.

\bibitem[H]{Heintze}
E.\ Heintze,
{\em Mannigfaltigkeiten negativer Kr\"ummung},
Ha\-bi\-li\-ta\-tions\-schrift, Universit\"at Bonn 1976.

\bibitem[Ja]{Jaco}
W.\ Jaco,
{\em Lectures on three-manifold topology},
Amer.\ Math.\ Soc.\ 43 (1980).

\bibitem[J-S]{JacoShalen}
W.\ Jaco and P.\ Shalen,
{\em Seifert fibred spaces in 3-manifolds},
Mem.\ Amer.\ Math.\ Soc.\ 220 (1979).

\bibitem[Jo]{Johannson}
K.\ Johannson,
{\em Homotopy equivalences of 3-manifolds with boundary},
Springer LNM 761 (1979).

\bibitem[K-L1]{action}
M.\ Kapovich and B.\ Leeb,
{\em Actions of discrete groups on nonpositively curved spaces},
preprint 1994.

\bibitem[K-L2]{graph}
M.\ Kapovich and B.\ Leeb,
{\em On quasi-isometries of graph-manifold groups},
preprint 1994.

\bibitem[L]{thesis}
B.\ Leeb, {\em 3-manifolds with(out) metrics of nonpositive curvature},
PhD Thesis, University of Maryland, 1992.

\bibitem[L-S]{deco}
B.\ Leeb and P.\ Scott,
{\em A geometric characteristic splitting in all dimensions},
in preparation.

\bibitem[S]{Scott}
P.\ Scott,
{\em The geometries of 3-manifolds},
Bull.\ London Math.\ Soc.\ 15 (1983), 401-487.

\bibitem[Th1]{Thurston1}
W.\ Thurston,
{\em The geometry and topology of 3-manifolds},
lecture notes, Princeton University.

\bibitem[Th2]{Thurston2}
W.\ Thurston,
{\em Hyperbolic structures on 3-manifolds, I,}
Ann.\ of Math.\ 124 (1986), 203-246.


\end{thebibliography}
\end{document}